\documentclass{article}

\usepackage{amsmath,amssymb,amsthm}

\newtheorem{theorem}{Theorem}[section]
\newtheorem{definition}[theorem]{Definition}
\newtheorem{proposition}[theorem]{Proposition}
\newtheorem{corollary}[theorem]{Corollary}
\newtheorem{lemma}[theorem]{Lemma}
\newtheorem{fact}[theorem]{Remark}
\newtheorem{exemplu}[theorem]{Example}

\newcommand{\bdfn}{\begin{definition}}
\newcommand{\edfn}{\end{definition}}
\newcommand{\bthm}{\begin{theorem}}
\newcommand{\ethm}{\end{theorem}}
\newcommand{\bprop}{\begin{proposition}}
\newcommand{\eprop}{\end{proposition}}
\newcommand{\bcor}{\begin{corollary}}
\newcommand{\ecor}{\end{corollary}}
\newcommand{\blem}{\begin{lemma}}
\newcommand{\elem}{\end{lemma}}
\newcommand{\bfact}{\begin{fact}}
\newcommand{\efact}{\end{fact}}
\newcommand{\bex}{\begin{exemplu}\begin{rm}}
\newcommand{\eex}{\end{rm}\end{exemplu}}

\def\N{{\mathbb N}}

\newcommand{\eps}{\varepsilon}

\newcommand{\bce}{\begin{center}}
\newcommand{\ece}{\end{center}}
\newcommand{\bi}{\begin{itemize}}
\newcommand{\ei}{\end{itemize}}
\newcommand{\be}{\begin{enumerate}}
\newcommand{\ee}{\end{enumerate}}
\newcommand{\bt}{\begin{tabular}}
\newcommand{\et}{\end{tabular}}
\newcommand{\beq}{\begin{equation}}
\newcommand{\eeq}{\end{equation}}
\newcommand{\ba}{\begin{array}} 
\newcommand{\ea}{\end{array}}
\newcommand {\bea} {\begin{eqnarray}}
\newcommand {\eea} {\end {eqnarray}}
\newcommand {\bua} {\begin{eqnarray*}}
\newcommand {\eua} {\end {eqnarray*}}

\newcommand{\ds}{\displaystyle}

\newcommand{\ra}{\rightarrow}

\begin{document}

\title{A quantitative Mean Ergodic Theorem for uniformly convex Banach spaces}

\author{ U. Kohlenbach${}^1$, L. Leu\c stean$^{1,2}$\\[0.2cm]
\footnotesize ${}^1$ Department of Mathematics, Technische Universit\" at Darmstadt,\\
\footnotesize Schlossgartenstrasse 7, 64289 Darmstadt, Germany\\[0.1cm]
\footnotesize${}^2$ Institute of Mathematics "Simion Stoilow'' of the Romanian Academy, \\
\footnotesize Calea Grivi\c tei 21, P.O. Box 1-462, Bucharest, Romania\\[0.1cm]
\footnotesize E-mails: kohlenbach,leustean@mathematik.tu-darmstadt.de
}
\maketitle

\begin{abstract}
\noindent We provide an explicit uniform bound on the local stability of 
ergodic averages in uniformly convex Banach spaces. Our result can also 
be viewed as a finitary version in the sense of T. Tao 
of the Mean Ergodic Theorem for such 
spaces and so generalizes similar results obtained for Hilbert 
spaces by Avigad, Gerhardy and  Towsner \cite{AviGerTow} and T. Tao \cite{Tao(07a)}.
\end{abstract}

\section{Introduction}

In the following $\N:=\{ 1,2,3,\ldots\}.$ \\ 
Let $X$ be a Banach space and $T:X\to X$. The {\em Cesaro mean} starting with $x\in X$ is the sequence $(x_n)_{n\geq 1}$ defined by $\ds x_n:=\frac{1}{n}\sum_{i=0}^{n-1}T^ix $.

In 1939, Garrett 
Birkhoff proved the following generalization of von Neumann's Mean Ergodic Theorem.

\bthm\cite{Bir39}\label{MET-Birkhoff}
Let $X$ be a uniformly convex Banach space and $T:X\to X$ be a linear operator 
with $\| Tx\| \le \| x\|$ for all $x\in X.$ 
Then for any $x\in X$, the Cesaro mean $(x_n)$ is convergent.
\ethm
 
In \cite{AviGerTow}, Avigad, Gerhardy and Towsner address the issue of finding an effective rate of convergence 
for $(x_n)$ in Hilbert spaces. They show that even for the separable 
Hilbert space $L_2$ there are simple computable such operators 
$T$ and computable points $x\in L_2$ such that there is no 
computable rate of convergence of $(x_n)$. In such a situation the best 
one can hope for is an effective bound on the following reformulation of 
the Cauchy property of $(x_n)$ which in logic is called the Herbrand normal 
form of the latter:  
\beq\label{Herbrand} 
\forall \varepsilon >0 \,\forall g:\N\to\N\,\exists n\in\N\, \forall 
i,j\in [n,n+g(n)] \ 
(\| x_i-x_j\| < \varepsilon). \eeq
It is trivial to see that (\ref{Herbrand}) is implied by the Cauchy 
property. However, ineffectively, also the converse implication holds. 
The mathematical relevance of this reformulation of convergence was 
recently pointed out by T. Tao (\cite{Tao(07),Tao(07a)}), who 
also uses the term `metastability'. 
In \cite{Koh05} (and refined in \cite{Gerhardy/Kohlenbach3}) 
a general logical metatheorem is 
proved that guarantees (given a proof of (\ref{Herbrand})) 
the extractability of an effective bound 
$\Phi(\eps,g,b)$ on `$\exists n$' in (\ref{Herbrand}) that is 
highly uniform in the sense that it only depends on $g,\varepsilon$ 
and an upper bound 
$\N\ni b\ge \| x\|$ but otherwise is independent from 
$x,X$ and $T.$ In fact, by a simple renorming argument one can 
always achieve to have the bound to depend on $b,\varepsilon$ only 
via $b/\varepsilon.$ The proof of this metatheorem, which is based 
on a recent extension and refinement of a technique from logic called G\"odel 
functional interpretation, provides an algorithm for extracting an explicit 
such $\Phi$ from a given proof (for a book treatment of all this 
see \cite{Kohlenbach(newcourse)}). Guided by this approach, 
Avigad, Gerhardy and Towsner \cite{AviGerTow} extract such a bound from a standard textbook proof of von Neumann's Mean Ergodic Theorem. 
A less direct proof for the existence of a bound with the above mentioned 
uniformity features is -- for a particular finitary dynamical system -- 
also given by T. Tao \cite{Tao(07a)}.   \\[2mm] 
In this note we apply the same methodology to Birkhoff's proof of  
theorem \ref{MET-Birkhoff} and extract an even easier to state bound for 
the more general case of uniformly convex Banach spaces. In this setting, 
the bound additionally depends on a given modulus of uniform convexity 
for $X$. Despite of our result being significantly more general 
then the Hilbert space case treated in \cite{AviGerTow}, the extraction of 
our bound is considerably more easy compared to \cite{AviGerTow} and 
even numerically better. 

\section{Main results}

Uniformly convex Banach spaces were introduced in 1936 by Clarkson in his seminal paper \cite{Cla36}. 

A Banach space $X$ is called {\em uniformly convex} if for all $\eps\in (0,2]$ there exists $\delta\in(0,1]$  such that for  all  $x,y\in X$, 
\beq
\| x\| \le 1, \quad \|y\|\le 1 \text{~~and~~}  \|x-y\|\geq \varepsilon  \text{~~imply~~} \left\|\frac12(x+y)\right\|\leq  1-\delta. \label{uc-def}
\eeq
A mapping $\eta:(0,2]\to (0,1]$ providing such a $\delta:=\eta(\eps)$ for given $\eps\in(0,2]$ is called a {\em modulus of uniform convexity}. 

Since the condition (\ref{uc-def}) is empty for $\eps>2$ we can simply 
extend any such $\eta$ to all strictly positive real numbers by 
stipulating $\eta'(\eps):=\eta(\min(2,\eps))$ if $\eta$ is not already 
defined for $\varepsilon >2.$ We will make free use of this without further 
mentioning. 
 
An example of a modulus of uniform convexity is Clarkson's {\em modulus of convexity} \cite{Cla36}, defined for any  Banach space $X$ as the function $\delta_X:[0,2]\to[0,1]$ given by
\beq
\delta_X(\varepsilon)=\inf\left\{1-\left\|\frac{x+y}{2}\right\|: \|x\|\le 1, \|y\|\le 1, \|x-y\|\ge \varepsilon \right\}.
\eeq
It is easy to see that $\delta_X(0)=0$ and that $\delta_X$ is nondecreasing. 
A well-known result is the fact that a Banach space $X$ is uniformly convex 
if and only if $\delta_X(\eps)>0$ for $\eps\in(0,2]$. Note that for uniformly convex spaces $X$, $\delta_X$ is the largest modulus of uniform convexity.

The  main result of our paper is a quantitative version of Birkhoff's generalization to uniformly convex Banach spaces of von Neumann's Mean Ergodic Theorem.

\bthm\label{quant-B-VN}
Assume that $X$ is a uniformly convex Banach space, $\eta$ is a modulus of uniform convexity  and  $T:X\to X$ is a linear operator 
with $\| Tx\| \le \| x\|$ for all $x\in X.$ 
Let $b>0$. Then for all $x\in X$ with $\|x\|\leq b$,
\beq
\forall\eps>0\,\forall g:\N\to\N\,\exists P\leq \Phi(\eps,g,b,\eta)\,\forall i,j\in[P,P+g(P)]\,\big(\|x_i-x_j\|<\eps\big).
\eeq
where $(x_n)$ is the Cesaro means starting with $x$ and 
\beq
\Phi(\eps,g,b,\eta):=M\cdot\tilde{h}^K(1),
\eeq
with
$$\ba{l}
\ds M:=\left\lceil\frac{16b}\eps\right\rceil, \gamma:=\frac\eps{16}
\eta\left(\frac{\eps}{8b}\right), \quad \ds K:=\left\lceil 
\frac{b}{\gamma}\right\rceil,\\[0.1cm]
h,\,\tilde{h}:\N\to\N, \,\,h(n):=2(Mn+g(Mn)),  \quad \ds 
\tilde{h}(n):=\max_{i\leq n}h(i).\\[0.1cm]
 \quad 
\ds 
\ea
$$
If $\eta(\varepsilon)$ can be written as $\ds \varepsilon\cdot 
\tilde{\eta}(\varepsilon)$ with $0<\varepsilon_1\le\varepsilon_2\to 
\tilde{\eta}(\varepsilon_1)
\le \tilde{\eta}(\varepsilon_2),$ then we can replace $\eta$ by 
$\tilde{\eta}$ and the constant `$16$' by `$8$' in the 
definition of $\gamma$  in the bound above.
\ethm

\bfact
Note that our bound $\Phi$ is independent from $T$ and 
depends on the space $X$ and the starting 
point $x\in X$ only via the modulus of convexity $\eta$ and 
the norm upper bound $b\ge \| x\|.$ Moreover, it is easy 
to see that the bound 
depends on $b$ and $\varepsilon$ only via $b/\varepsilon$.
\efact

As an immediate consequence of our theorem we get a quantitative version of 
von Neumann's Mean Ergodic Theorem.

\bcor\label{quant-Hilbert}
Assume that $X$ is a Hilbert space and  $T:X\to X$ is a linear operator 
with $\| Tx\| \le \| x\|$ for all $x\in X.$ 
Let $b>0$. Then for all $x\in X$ with $\|x\|\leq b$,
\beq
\forall\eps>0\,\forall g:\N\to\N\,\exists P\leq 
\Phi(\eps,g,b)\,\forall i,j\in[P,P+g(P)]\,\big(\|x_i-x_j\|<\eps\big).
\eeq
where $(x_n),$ $\Phi$ are defined as above, but with 
$\ds K:=\left\lceil \frac{512b^2}{\eps^2}\right\rceil$. 
\ecor
\begin{proof}
It is well-known that as a modulus of uniform convexity of a 
Hilbert space $X$ we can take $\eta(\varepsilon):=\eps^2/8$ with
$\tilde{\eta}(\eps):=\eps/8$ satisfying the requirements in the last claim 
of the theorem. 
\end{proof}
We get a similar result for $L_p$-spaces ($2<p<\infty$), using the fact 
that $\ds \eta(\eps)=\frac{\varepsilon^p}{p2^p}$ is a modulus of uniform convexity for $L_p$ 
(see e.g. \cite{Kohlenbach(asymptotic)}). 
Note that  $\ds \frac{\varepsilon^p}{p2^p} =\varepsilon \cdot
\tilde{\eta}_p(\varepsilon)$ with $\ds \tilde{\eta}_p(\varepsilon)=\frac{\varepsilon^{p-1}}{p2^p} $
satisfying the monotonicity condition in the theorem above.
\bfact The bound extracted in \cite{AviGerTow} for Hilbert spaces  
is the following one 
\[ \begin{array}{l} \Phi(\eps,g,b)= h^K(1)\ \mbox{where} \ 
h(n)= n+2^{13}
\rho^4\tilde{g}((n+1)\tilde{g}(2n\rho)\rho^2) \ 
\mbox{with} \\  
\rho=\left\lceil \frac{b}{\varepsilon}\right\rceil, \ K=2^9\rho^2, \ 
\mbox{and} \ \ds \tilde{g}(n)=\max_{i\le n}(i+g(i)). \end{array} \]
Note that the number of iterations essentially is the same as in our bound 
in 
corollary \ref{quant-Hilbert} above but that the function being iterated in 
our corollary is much simpler. Roughly speaking, our bound for the 
general nonexpansive case (i.e. 
$\| Tx\| \le \| x\|$ for all $x\in X$) corresponds to 
the one obtained in \cite{AviGerTow} for the special case of $T$ being an 
isometry with $\Phi$ as above but 
\[ h(n)= n+2^{13}
\rho^4\tilde{g}((n+1)\tilde{g}(1)\rho^2) \] 
from which in \cite{AviGerTow} it is derived that 
$\Phi(\eps,g,b)=2^{O(\rho^2\log \rho)}$ (with $\rho:=\lceil \frac{b}{\eps}\rceil$) for linear functions $g,$ i.e. $g=O(n).$ \\ Our corollary 
\ref{quant-Hilbert} generalizes this to $T$ being nonexpansive 
rather than being an isometry.
  
\efact

\section{Technical lemmas}

\blem\label{GLB-lemma}
Let  $(a_n)_{n\geq 0}$ be a sequence of nonnegative real numbers. Then
\begin{itemize}
\item[(i)] $\ds \forall\eps>0\, \forall g:\N\to\N\,\exists 
N\leq\Theta(b,\eps,g)\,\big(a_N\leq a_{g(N)}+ \eps\big), $\\
where $\ds \Theta(b,\eps,g):= \max_{i\leq K}g^i(1), \,\, 
b \geq  a_0,  \,\, K :=\left\lceil \frac{b}{\varepsilon} 
\right\rceil$. Moreover, $N=g^i(1)$ 
for some $i<K$.
\item[(ii)] $\ds \forall\eps>0\, 
\forall g:\N\to\N\,\exists N\leq h^K(1)\,\forall 
m\leq g(N)\big(a_N\leq a_m+\eps\big),$\\
where $\ds  h(n):=\max_{i\leq n}g(i)$ and  $b, K$ are as above.
\end{itemize}
\elem
\begin{proof}
\begin{itemize}
\item[(i)]  See e.g. \cite[Lemma 6.3]{KohLeu08}
\item[(ii)]  Let $\eps>0$, $g:\N\to \N$ and define
\bua
\tilde{g}:\N\to\N, \quad \tilde{g}(n):= 
\text{~the least~} i\leq g(n) \text{~satisfying~} 
a_i=\min\{a_j\mid j\leq g(n)\}. 
\eua
Then, for all $n\in\N$ and for all $m\leq g(n)$, we have that  
$a_m\geq a_{\tilde{g}(n)}$. Applying now (i) for $\eps$ and 
$\tilde{g}$, we get that there exists $N\leq \Theta(b,\eps,\tilde{g})$ 
such that $a_N\leq a_{\tilde{g}(N)}+ \eps\leq a_m +\eps$ for all 
$m\leq g(N)$. Let us now define $\ds h:\N\to\N, \,\,
h(n)=\max_{i\leq n}g(i)$. Then $h$ is nondecreasing and 
$h(n)\geq g(n)\geq \tilde{g}(n)$ for all $n\in\N$. It is easy to see 
$h^i(n)\geq \tilde{g}^i(n)$ and $h^{i}(1)\ge h^{i-1}(1)$ for all 
$i,n\in\N.$ Hence, 
$\ds h^K(1)=\max_{i\leq K}h^i(1)\geq 
\max_{i\leq K}\tilde{g}^i(1)=\Theta(b,\eps,\tilde{g})\geq N$.
\end{itemize}
\end{proof}

\blem\label{uc-u}
Let $X$ be a uniformly convex Banach space and $\eta$ be a modulus of uniform convexity. Define $u_\eta:(0,2]\ra(0,1], \quad \ds u_\eta(\eps)=\frac{\eps}2\cdot\eta\left(\eps\right)$. Then for all  $\eps>0$  and for all $x,y\in X$
\beq
\|x\|\leq\|y\|\leq 1 \text{~and~} \|x-y\| \geq \eps \quad \text{imply}\quad \left\|\frac{1}{2}(x+y)\right\|\leq \|y\|-u_\eta(\eps). 
\eeq 
We use the notation $u_X$ for $\ds u_{\delta_X}$, where $\delta_X$ is the modulus of convexity. \\ 
If $\eta(\varepsilon)$ can be written as $\varepsilon\cdot \tilde{\eta}(\varepsilon)$ with $0<\varepsilon_1\le\varepsilon_2\to \tilde{\eta}(\varepsilon_1)
\le \tilde{\eta}(\varepsilon_2),$ then we can replace $u_{\eta}$ by 
$\tilde{u}_{\eta}(\eps):=\eps\cdot\tilde{\eta}(\eps).$
\elem
\begin{proof}
We have that $\ds\frac{\|x\|}{\|y\|}\leq \frac{\|y\|}{\|y\|}=1$ and $\ds\frac1{\|y\|}\|x-y\|\geq \frac{\eps}{\|y\|}\geq\eps$, since $\|y\|\leq 1$. Applying the fact that $\eta$ is a modulus of uniform convexity, we get that $\ds \frac1{\|y\|}\left\|\frac{1}{2}(x+y)\right\|\leq 1 -\eta(\eps)$, hence
\bua
\left\|\frac{1}{2}(x+y)\right\|\leq \|y\|-\|y\|\eta(\eps)\leq \| y\|- 
u_\eta(\eps),
\eua
since $\ds \|y\|\geq \frac12(\|x\|+\|y\|)\geq \frac12\|x-y\|\geq \frac{\eps}2$.
\\ The last claim follows from 
\bua 
\left\|\frac{1}{2}(x+y)\right\|\leq \|y\|-\|y\|\eta(\eps/\| y\|)
=\| y\| -\eps\cdot \tilde{\eta} (\eps/\| y\|)\leq \| y\| 
-\eps\cdot \tilde{\eta}(\eps).\eua
\end{proof}

The following lemma collects some facts already remarked by Birkhoff in his paper \cite{Bir39}. For completeness, we give the proofs here.
\blem\cite{Bir39}\label{linear-Cesaro}
Let $X$ be a Banach space, $T:X\to X$ be linear and $(x_n)$ be the Cesaro mean starting with $x$. 
\begin{itemize}
\item[(i)]  For all $n,k\in\N$,
\bea
x_{n+k}&=& \frac{n}{n+k}x_n+\frac{1}{n+k}\sum_{i=0}^{k-1}T^{n+i}x,\label{x-n+k=}\\
x_{kn}&=& \frac{1}{k}\sum_{i=0}^{k-1}T^{in}x_n,\label{linear-cesaro-2}\\
x_{2kn}&=& \frac{1}{k}\sum_{i=0}^{k-1}\frac{1}{2}T^{in}\left(x_n+T^{kn}x_n\right). \label{linear-cesaro-1}
\eea
\item[(ii)]  Assume moreover that $T$ satisfies $\| Tx\| \le \| x\|$ 
for all $x\in X.$ Then for all  $n,k\in\N$,
\bea
\|x_{n+k}-x_n\|&\leq & \frac{2k\|x\|}{n+k}\label{x-n+k-x-n},\\
\|x_{kn}-x_n\|&\leq & \max_{i=0,\ldots,k-1}\|T^{in}x_n-x_n\| \label{x-kn-x-n}
\eea
\end{itemize}
\elem
\begin{proof}
\begin{itemize}
\item[(i)] (\ref{x-n+k=}) is obvious, (\ref{linear-cesaro-2}) and (\ref{linear-cesaro-1}) are obtained by grouping terms:
\bua
x_{kn}&=& \frac1{kn}\sum_{j=0}^{nk-1}T^jx=\frac1{kn}\sum_{i=0}^{k-1}\left(T^{in}x+T^{in+1}x+\ldots + T^{in+(n-1)}x\right)\\
&=&\frac1k\sum_{i=0}^{k-1}T^{in}\left(\frac1n(x+\ldots + T^{n-1}x)\right)=\frac{1}{k}\sum_{i=0}^{k-1}T^{in}x_n\\[0.1cm]
x_{2kn}&=& \frac1{2kn}\sum_{j=0}^{2nk-1}T^jx=\frac1{2kn}\sum_{i=0}^{k-1}\left(\sum_{j=0}^{n-1}T^{in+j}x+\sum_{j=0}^{n-1}T^{(k+i)n+j}x\right)\\
&=&\frac1k\sum_{i=0}^{k-1}\frac12T^{in}\left(\frac1n\sum_{j=0}^{n-1}T^jx+T^{kn}\left(\frac1n\sum_{j=0}^{n-1}T^{j}x\right)\right)\\
&=&\frac{1}{k}\sum_{i=0}^{k-1}\frac{1}{2}T^{in}\left(x_n+T^{kn}x_n\right).
\eua
\item[(ii)] By assumption we 
have that $\|Ty\|\leq \|y\|$ for all $y\in X$, so  $\|T^nx\|\leq \|x\|$ for all $n\in\N$ and, moreover, $\ds\|x_n\|\leq\frac1n\sum_{i=0}^{n-1}\|T^ix\|\leq \|x\|$ for all $n\in\N$.
\bua
\|x_{n+k}-x_n\|&\stackrel{(\ref{x-n+k=})}{=}& \left\|\left(\frac{n}{n+k}x_n+\frac{1}{n+k}\sum_{i=0}^{k-1}T^{n+i}x\right)-x_n\right\|\\
&=& \left\|\frac{-k}{n+k}x_n+\frac{1}{n+k}\sum_{i=0}^{k-1}T^{n+i}x\right\|\\
&\leq& \frac{2k}{n+k}\|x\|.
\eua
\bua
\|x_{kn}-x_n\|&\stackrel{(\ref{linear-cesaro-2})}{=}& \left\|\frac{1}{k}\sum_{i=0}^{k-1}T^{in}x_n-x_n\right\|=\left\|\frac{1}{k}\sum_{i=0}^{k-1}(T^{in}x_n-x_n)\right\|\\
&\leq &\frac{1}{k}\sum_{i=0}^{k-1}\|T^{in}x_n-x_n\|\leq 
\max_{i=0,\ldots,k-1}\|T^{in}x_n-x_n\|.
\eua
\end{itemize}
\end{proof}

\section{Proof of Theorem \ref{quant-B-VN}}

Let $x\in X$, $\eps>0$ and $g:\N\to\N$ be arbitrary and  
$\Phi,b,M,\gamma,K,h,\tilde{h}$ as in the hypotheses. 
Then $\ds M\geq \frac{16b}\eps$, that is $\ds \frac{2b}M\leq\frac\eps{8}$. 

Let $N$ be obtained by applying Lemma \ref{GLB-lemma}(ii) for the sequence 
$(\|x_n\|)_{n\geq 1}$ and  the above $\gamma$ and $h$. It follows 
that $0<N\leq \tilde{h}^K(1)$ exists satisfying
\beq
\forall m\leq h(N)(\|x_N\|\leq \|x_m\|+\gamma).\label{ND-GLB-ar}
\eeq
 Denote for all $k\in\N$,
\beq
y_k:=\|T^{kN}x_N-x_N\|.
\eeq
{\bf Claim:
 For all $\ds k\leq \frac{h(N)}{2N}$, we have that $\ds y_k \leq 
\frac\eps{8}.$}\\[0.1cm]
{\bf Proof of claim:} If $y_k=0$, then it is obvious, so we can assume 
in the sequel that $y_k\neq 0$. We get that for all $k\in\N$ 
\bua
\left\|\frac1b T^{kN}x_N\right\|\leq \left\|\frac1b x_N\right\|\leq 
\frac{\|x\|}b\leq 1 \text{~~and~~}\\  
\frac{y_k}{b}=\left\|\frac1b(T^{kN}x_N-x_N)\right\| 
\leq \frac1b(\|T^{kN}x_N\|+\|x_N\|)\leq2\frac{\|x\|}{b}\leq 2.
\eua
Thus, applying Lemma \ref{uc-u}, we get that
\beq
\left\|\frac{1}{2b}(T^{kN}x_N+x_N)\right\| \leq 
\frac1b\|x_N\|- u_X\left(\frac{y_k}{b}\right),
\eeq
that is
\beq
\left\|\frac12(T^{kN}x_N+x_N)\right\|\leq 
\|x_N\|-bu_X\left(\frac{y_k}{b}\right)
\eeq
for all $k\in\N$.

Using now (\ref{linear-cesaro-1}) of Lemma \ref{linear-Cesaro}, we obtain 
\bua
\|x_{2kN}\| &=&\left\|\frac{1}{k}
\sum_{i=0}^{k-1}\frac{1}{2}T^{iN}
\left(x_N+T^{kN}x_N\right)\right\| 
\leq \frac{1}{k}\sum_{i=0}^{k-1}
\left\|T^{iN}\left(\frac{1}{2}(x_N+T^{kN}x_N)\right)\right\|\\
&\leq & \frac{1}{k}\sum_{i=0}^{k-1}\left\|\frac{1}{2}(x_N+T^{kN}x_N)
\right\|=\left\|\frac{1}{2}(x_N+T^{kN}x_N)\right\|\\
&\leq & \|x_N\|-bu_X\left(\frac{y_k}{b}\right).
\eua
On the other hand, applying (\ref{ND-GLB-ar}), we get for 
$\ds k\leq \frac{h(N)}{2N}$\bua
\|x_{2kN}\| \geq \|x_N\|-\gamma.
\eua
Thus we must have that 
\beq
bu_X\left(\frac{y_k}{b}\right)\leq\gamma \quad 
\text{for all~~}k\leq\frac{h(N)}{2N}. \label{yk-gamma}
\eeq
Assume that $\ds y_k> \frac\eps{8}$. Then, since 
$\delta_X$ is nondecreasing and $\delta_X\geq \eta$, we get that
\beq \label{tilde-line}
bu_X\left(\frac{y_k}{b}\right)= 
b\cdot\frac{y_k}{2b}\cdot\delta_X\left(\frac{y_k}{b}\right)> 
\frac{\eps}{16}\delta_X\left(\frac{\eps}{8b}\right)\geq 
\frac{\eps}{16}\eta\left(\frac{\eps}{8b}\right) =\gamma,
\eeq
that is a contradiction with (\ref{yk-gamma}). Hence, we must have 
$\ds y_k\leq \frac\eps{8}$ for all $\ds k\leq\frac{h(N)}{2N}$.
This finishs the proof of the claim. \\[1mm]
Using the claim it 
follows that for all $\ds 0< m\leq \frac{h(N)}{2N}$ and $0\leq i<N$, 
we get that
\beq
\|x_{mN+i}-x_N\| \leq  \frac{2b}m+\frac{\eps}8, \label{ineq-m-i}
\eeq
since
\bua
\|x_{mN+i}-x_N\|&\leq & \|x_{mN+i}-x_{mN}\|+\|x_{mN}-x_N\|\\
&\leq & \frac{2ib}{mN+i}+\|x_{mN}-x_N\|, \quad \text{by (\ref{x-n+k-x-n}) 
and the fact that~~} \|x\|\leq b\\
&< & \frac{2b}m+\|x_{mN}-x_N\|, \quad \text{since~~}0\leq i<N 
\text{~~implies} \frac{2i}{mN+i} <\frac2m\\
&\leq & \frac{2b}m+\max_{{j=0,\ldots,m-1}}y_j, \quad 
\text{by (\ref{x-kn-x-n})}\\
&\leq & \frac{2b}m+\frac{\eps}8 \quad \text{by the above claim.}
\eua
Let us define $P:=MN\leq \Phi(\eps,g,b,\eta)$ and take $j\in[P,P+g(P)]$. 
Then  there are $q\in\N_0, 0\leq i <N$ such that $j-P=Nq+i$; 
moreover $Nq\leq j-P\leq g(P)=g(MN)$, so $\ds q\leq\frac{g(MN)}{N}$.
It follows that
\bua
\|x_j-x_P\|&=& \|x_{MN+Nq+i}-x_{MN}\| =\|x_{N(M+q)+i}-x_{MN}\|\\
&\leq & \|x_{N(M+q)+i}-x_{N}\|+\|x_{MN}-x_N\| \\
& < & \frac{2b}{M+q}+\frac{\eps}8+\frac{2b}M+\frac{\eps}8\leq 
\frac{\eps}4+\frac{4b}M\leq \frac{\eps}2,
\eua
since $M\leq M+q\leq M+\frac{g(MN)}{N}=\ds \frac{h(N)}{2N}$, so we can 
apply (\ref{ineq-m-i}) with $m:=M$ and $m:=M+q$.

It follows immediately that for all $j,l\in[P,P+g(P)]$, we have that
\[
\|x_j-x_l\|\leq \|x_j-x_P\|+\|x_l-x_P\|<\eps.
\]
The last claim of the theorem follows using the last claim in Lemma 
\ref{uc-u} with 
$\ds \gamma:=\frac{\eps}{8}\tilde{\eta}\left(\frac{\eps}{8b}\right)$ 
and  
$\ds \tilde{u}_{\eta}$ instead of $u_X.$ Then (\ref{tilde-line}) 
needs to be replaced by 
\bua
b\cdot 
\tilde{u}_{\eta}\left(\frac{y_k}{b}\right)= y_k\cdot\tilde{\eta}
\left(\frac{y_k}{b}\right)> \frac{\eps}{8}\tilde{\eta}
\left(\frac{\eps}{8b}\right)= \gamma.
\eua
{\bf Final remark on the extraction of the bound:} The only ineffective 
principle used in Birkhoff's original proof is the fact that any 
sequence $(a_n)$ of positive real numbers has an infimum (GLB).
In our analysis we first replaced this analytical existential statement 
by a purely arithmetical one, namely 
\[ \mbox{(GLB$_{ar}$)}: \ 
\forall \eps>0\,\exists n\in\N\,\forall m\in\N\ (a_n\le a_m+\eps). \] 
This principle still is ineffective as there (in general) is no computable 
bound on `$\exists n\in\N$' (even for computable $(a_n)).$ We then 
carried out (informally) a version of G\"odel's functional interpretation  
by which (GLB$_{ar}$) gets replaced in the proof by the quantitative form 
provided in lemma \ref{GLB-lemma}. 
For the general underlying facts from logic that 
guarantee this to be possible see \cite{Kohlenbach(newcourse)}.

\end{document}